\documentstyle[12pt]{article}
\def\Box{{\setlength{\unitlength}{1.1 ex}
\begin{picture}(1,1)(-0.2,-0.2)
\put(0,0){\framebox(1,1){}}
\end{picture}}}

\font\elevenbb=msbm10 at 10.95pt

\def\N{\hbox{\elevenbb N}}
\def\M{\hbox{\elevenbb M}}
\def\R{\hbox{\elevenbb R}}

\def\Gr{Gr\"obner }
\def \bg #1 {\begin{tabular}{{#1}}}
\def \nd {\end{tabular}}
\newcommand \mwhile {{\bf while}\hspace{0.3cm}}
\newcommand \mrepeat {{\bf repeat}\hspace{0.3cm}}
\newcommand \muntil {{\bf until}\hspace{0.3cm}}
\newcommand \mfore {{\bf for\hspace{0.3cm}each}\hspace{0.3cm}}
\newcommand \mdo {{\bf do}\hspace{0.3cm}}

\newcommand \mif {{\bf if}\hspace{0.3cm}}
\newcommand \mthen {{\bf then}\hspace{0.3cm}}
\newcommand \melse {{\bf else}\hspace{0.3cm}}
\newcommand \mchoose {{\bf choose}\hspace{0.3cm}}
\newcommand \mand {{\bf and}\hspace{0.3cm}}

\newcommand \mbegin {{\bf begin}}
\newcommand \mend {{\bf end}}
\newcommand \bb {\hspace{0.3cm}}
\newcommand \h {\hspace{0.5cm}}
\newcommand \hh {\hspace{1.0cm}}
\newcommand \hhh {\hspace{1.5cm}}
\newcommand \hhhh {\hspace{2.0cm}}
\newcommand \hhhhh {\hspace{2.5cm}}
\newcommand \hhhhhh {\hspace{3.0cm}}
\newcommand \hhhhhhh {\hspace{3.5cm}}
\newcommand \hhhhhhhh {\hspace{4.0cm}}

\newcounter{cc}
\newcommand \hln {\hfill \addtocounter{cc}{1} \arabic{cc}
\vskip 0.0cm \noindent }
\newtheorem{definition}{Definition}[section]
\newtheorem{corollary}[definition]{Corollary}
\newtheorem{remark}[definition]{Remark}
\newtheorem{proposition}[definition]{Proposition}
\newtheorem{example}[definition]{Example}
\newtheorem{theorem}[definition]{Theorem}

\begin{document}
\title{\bf Involutive Division Technique: Some Generalizations and
 Optimizations}
\author{Vladimir P. Gerdt \\
       Laboratory of Computing Techniques and Automation\\
       Joint Institute for Nuclear Research\\
       141980 Dubna, Russia \\
       gerdt@jinr.ru}
\date{}
\maketitle
\begin{abstract}
In this paper, in addition to the earlier introduced involutive
divisions, we consider a new class of divisions induced by
admissible monomial orderings. We prove that these divisions are
noetherian and constructive. Thereby each of them allows one to
compute an involutive \Gr basis of a polynomial ideal by
sequentially examining multiplicative reductions of
nonmultiplicative prolongations. We study dependence
of involutive algorithms on the completion ordering. Based on
properties of particular involutive divisions two computational
optimizations are suggested. One of them consists in a special
choice of the completion ordering. Another optimization is related
to recomputing multiplicative and nonmultiplicative variables in
the course of the algorithm.
\end{abstract}

\section{Introduction}
\noindent
In paper~\cite{GB1} a concept of involutive monomial division was
invented which forms the foundation of general involutive
algorithms~\cite{GB1,GB2} for construction of \Gr
bases~\cite{Buch65} of a special form called involutive. This
notion, by a well-known correspondence~\cite{Pommaret94,Gerdt95}
between polynomials and linear homogeneous partial differential
equations (PDEs) with constant coefficients, follows the notion of
involutivity for PDEs. An involutive form of a system of PDEs is
its interreduced completion by the differential consequences called
prolongations\footnote{Prolongation for PDE means its
differentiation whereas for a polynomial this means multiplication
by the corresponding variable.}, incorporating all integrability
conditions into the system~\cite{Janet, Pommaret78}. The
integrability conditions play the same role in the completion
procedure for PDEs as nontrivial $S$-polynomials in the Buchberger
algorithm~\cite{Buch85,BWK93} for construction of \Gr bases.

Given a finite polynomial set and an admissible monomial ordering, an
involutive division satisfying the axiomatic properties proposed
in~\cite{GB1} leads to a self-consistent separation of variables for
any polynomial in the set into disjoint subsets of so-called
multiplicative and nonmultiplicative variables.

The idea of the separation of variables into multiplicative and
nonmultiplicative goes back to classical papers of
Janet~\cite{Janet} and Thomas~\cite{Thomas}. They used particular
separations of independent variables for completing systems of
partial differential equations to involution. More recently one of
the possible separations already introduced by Janet~\cite{Janet}
was intensively used by Pommaret~\cite{Pommaret78} in the formal
theory of partial differential equations. These classical
separations allow one to generate the integrability conditions by
means of multiplicative reductions of nonmultiplicative
prolongations. Just this fact was first used in~\cite{ZB93} as a
platform for an involutive algorithm for construction of Pommaret
bases of polynomial ideals.

If an involutive division satisfies some extra conditions:
noetherity and constructivity~\cite{GB1}, then an involutive basis
may be constructed algorithmically by sequential examination of
single nonmultiplicative prolongations only. Whereas Thomas and
Janet divisions satisfy all the extra conditions, Pommaret
division, being constructive, is non-noetherian. This implies that
Pommaret bases of positive dimensional ideals may be infinite. The
uniqueness properties of involutive bases are investigated
in~\cite{GB2} where a special form of an algorithm proposed for
construction of a minimal involutive basis which is unique much
like to a reduced \Gr basis. In addition to the above mentioned
classical divisions, in paper~\cite{GB2} two more divisions were
introduced which satisfy all the extra conditions.

Recently it was shown~\cite{Apel} that one can also construct
different possible separations of variables for a fixed monomial
set. These separations can not be considered, generally, as
functions of a set and its element as defined in~\cite{GB1}.
Nevertheless, the results of paper~\cite{Apel} demonstrate for a
wide class of divisions how one can change the division dynamically
in the course of the completion. This increases the flexibility of
the involutive technique and may also increase the efficiency of
computations.

An involutive basis is a special kind of \Gr basis, though,
generally, it may be redundant. However, extra elements in the
former may facilitate many underlying computations. The structure
of a Pommaret basis, for example, reveals a number of attractive
features convenient for solving zero-dimensional polynomial
systems~\cite{Z96}. An involutive basis for any division allows one
to compute easily the Hilbert function and the Hilbert polynomial
by explicit and compact formulae~\cite{Apel,GBC98}.

Computation of Janet bases relying upon the original Janet
algorithm was implemented in Reduce and used for finding the size
of a Lie symmetry group for PDEs~\cite{Schwarz} and for
classification of ordinary differential equations admitting
nontrivial Lie symmetries~\cite{Schwarz1}. The study of algorithmic
aspects of the general completion procedure for Pommaret division
and implementation in Axiom was done in~\cite{Seiler}.  The
completion to involution of polynomial bases for Pommaret division
was algorithmized and implemented in Reduce, first, in~\cite{ZB93},
and then with algorithmic improvements in~\cite{GB1}. The main
improvement is incorporation of an involutive analogue of
Buchberger's chain criterion. Recently different involutive
divisions were implemented also in Mathematica~\cite{GBC98}.

In the present paper we introduce a class of involutive divisions
induced by admissible orderings and prove their noetherity and
constructivity. For the new class of divisions, along with the
classical ones and two divisions of paper~\cite{GB2}, we study the
stability of the partial involutivity for monomial and polynomial
sets under their completion by irreducible nonmultiplicative
prolongations. We generalize the involutive algorithms to
different main and completion orderings. The completion ordering
serves for selection of a nonmultiplicative prolongation to be
treated next. In so doing, a completion ordering defines the
selection strategy in involutive algorithms similar to the
selection strategy for critical pairs in Buchberger
algorithm~\cite{Buch85,BWK93}. For different divisions we find some
completion orderings which preserve the property of partial
involutivity and thereby may save computing time. We indicate also a
'pairwise' property which is valid for all known divisions. This
property can be used to efficiently recompute the separations when
a new polynomial has to be added.

\section{Background of Involutive Method}

In this section, we recall basic definitions and facts of
papers~\cite{GB1,GB2} which are used in the next sections.

\subsection{Preliminaries}

Let $\N$ be the set of nonnegative integers, and $\M=\{x_1^{d_1}\cdots
x_n^{d_n}\ |\ d_i\in \N\}$ be a set of monomials in the polynomial ring
$\R=K[x_1,\ldots,x_n]$ over a field $K$ of characteristic zero.

By $deg(u)$ and $deg_i(u)$ we denote the total degree of $u\in \M$
and the degree of variable $x_i$ in $u$, respectively. For the
least common multiple of two monomials $u,v\in \M$ we shall use the
conventional notation $lcm(u,v)$. If monomial $u$ divides monomial
$v$ we shall write $u|v$. In this paper we shall distinguish two
admissible monomial orderings: {\em main ordering} and {\em
completion ordering} denoted by $\succ$ and $\sqsupset$,
respectively. The main ordering serves, as usually, for isolation
of the leading terms in polynomials whereas the completion ordering
is used for taking the lowest nonmultiplicative prolongations by
the normal strategy~\cite{GB1} and thereby controlling the property
of partial involutivity. Besides, throughout the paper we shall
assume that
\begin{equation}
x_1\succ x_2\succ\cdots\succ x_n\,. \label{var_order}
\end{equation}
The leading monomial and the leading coefficient of the polynomial $f\in
\R$ with respect to $\succ$ will be denoted by $lm(f)$ and $lc(f)$,
respectively. If $F\subset \R$ is a polynomial set, then by $lm(F)$ we
denote the leading monomial set for $F$, and $Id(F)$ will denote the ideal
in $R$ generated by $F$. The least common multiple of the set $\{lm(f)\ |\
f\in F\}$ will be denoted by $lcm(F)$.

\subsection{Involutive Monomial Division}

\begin{definition}
{\em An {\em involutive division} $L$ on $\M$ is given, if for any finite
monomial set $U\subset \M$ and for any $u\in U$ there is given a
submonoid $L(u,U)$ of $\M$ satisfying the conditions:
\begin{tabbing}
~~(a)~~\=If $w\in L(u,U)$ and $v|w$, then $v\in L(u,U)$. \\
~~(b)  \>If $u,v\in U$ and $uL(u,U)\cap vL(v,U)\not=\emptyset$, then
$u\in vL(v,U)$ or $v\in uL(u,U)$. \\
~~(c)  \> If $v\in U$ and $v\in uL(u,U)$, then $L(v,U)\subseteq L(u,U)$. \\
~~(d)  \> If $V\subseteq U$, then $L(u,U)\subseteq L(u,V)$ for all $u\in V$.
\end{tabbing}
Elements of $L(u,U)$, $u\in U$ are called {\em multiplicative} for $u$.  If
$w\in uL(u,U)$ we shall write $u|_L w$ and call $u$ {\em ($L-$)involutive
divisor} of $w$. The monomial $w$ is called {\em ($L-$)involutive multiple}
of $u$. In such an event the monomial $v=w/u$ is {\em multiplicative} for
$u$ and the equality $w=uv$ will be written as $w=u\times v$. If $u$ is a
conventional divisor of $w$ but not an involutive one we shall write, as
usual, $w=u\cdot v$. Then $v$ is said to be {\em nonmultiplicative} for
$u$. } \label{inv_div}
\end{definition}

\begin{definition}{\em We shall say that an involutive division $L$ is
 {\em globally defined} if for any $u\in \M$ its multiplicative monomials
 are defined irrespective of the monomial set $U\ni u$, that is, if
 $L(u,U)=L(u)$.
} \label{df_global}
\end{definition}

\noindent
Definition~\ref{inv_div} for every $u\in U$ provides the separation
\begin{equation}
\{x_1,\ldots,x_n\}=M_L(u,U)\cup NM_L(u,U),\quad M_L(u,U)\cap NM_L(u,U)=
\emptyset
\label{part}
\end{equation}
of the set of variables into two subsets: {\em multiplicative}
$M_L(u,U)\subset L(u,U)$ and {\em non\-mul\-ti\-pli\-ca\-ti\-ve}
$NM_L(u,U)\cap L(u,U)=\emptyset$.
Conversely, if for any finite set $U\subset \M$ and any $u\in U$ the
separation~(\ref{part}) is given such that the corresponding
submonoid $L(u,U)$ of monomials in variables in $M_L(u,U)$ satisfies
the conditions (b)-(d), then the partition generates an involutive
division.
The conventional monomial division, obviously, satisfies condition
(b) only in the univariate case.

In what follows monomial sets are assumed to be finite.

\begin{definition}{\em A monomial set $U\in \M$ is {\em involutively
 autoreduced} or {\em $L-$autoreduced} if the condition
 $uL(u,U)\cap vL(v,U)=\emptyset $ holds for all distinct $u,v\in U$.
}\end{definition}

\begin{definition}{\em Given an involutive division $L$, a
 monomial set $U$
 is {\em involutive} with respect to $L$ or $L-$involutive if
$$
 (\forall u\in U)\ (\forall w\in \M)\
 (\exists v\in U)\ \ [\ uw\in vL(v,U)\ ]\,.
$$
} \label{inv_mset}
\end{definition}

\begin{definition}{\em We shall call the set $\cup_{u\in U}\,u\,\M$ {\em the
cone} generated by $U$ and denote it by $C(U)$. The set $\cup_{u\in
U}\,u\,L(u,U)$ will be called {\em the involutive cone} of $U$ with
respect to $L$ and denoted by $C_L(U)$. } \label{cone}
\end{definition}
Thus, the set $U$ is $L-$involutive if its cone $C(U)$ coincides with
its involutive cone $C_L(U)$.

\begin{definition}{\em An $L-$involutive monomial set $\tilde{U}$
 is called $L-${\em completion} of a set $U\subseteq \tilde{U}$ if
$$
 (\forall u\in U)\ (\forall w\in \M)\
 (\exists v\in \tilde{U})\ \ [\ uw\in vL(v,\tilde{U})\ ]\,.
$$
 If there exists a finite $L-$completion $\tilde{U}$ of a finite set
 $U$, then the latter is {\em finitely generated} with respect
 to $L$. The involutive division $L$ is {\em noetherian} if
 every finite set $U$ is finitely generated with respect to $L$.
} \label{id_noetherian}
\end{definition}

\begin{proposition}\cite{GB1}
If an involutive division $L$ is noetherian, then every monomial ideal
has a finite involutive basis $\bar{U}$.
\label{pr_fin_ib}
\end{proposition}

\begin{definition}
 {\em A monomial set $U$ is called {\em locally involutive} with respect
 to the involutive division $L$ if
$$
(\forall u\in U)\ (\forall x_i\in NM_L(u,U))\ (\exists v\in U)\
\ [\ v|_L(u\cdot x_i)\ ]\,.
$$
}
\label{mon_loc_inv}
\end{definition}

\begin{definition}
{\em A division $L$ is called {\em continuous} if for any set $U\in
\M$ and for any finite sequence $\{u_i\}_{(1\leq i\leq k)}$ of
elements in $U$ such that
\begin{equation}
(\forall \,i< k)\ (\exists x_j\in NM_L(u_i,U))\ \ [\ u_{i+1}|_L u_i\cdot
x_j\ ] \label{cont_cond}
\end{equation}
the inequality $u_i\neq u_j$ for $i\neq j$ holds.
} \label{def_cont}
\end{definition}

\begin{theorem}\cite{GB1}
If an involutive division $L$ is continuous then local involutivity
of a monomial set $U$ implies its involutivity.
\label{th_cont}
\end{theorem}

\begin{definition}
{\em  A continuous involutive division $L$
 is {\em constructive} if for any $U\subset \M$, $u\in U$,
 $x_i\in NM_L(u,U)$ such that $u\cdot x_i$ has no involutive divisors
 in $U$ and
$$
 (\forall v\in U)\ (\forall x_j\in NM_L(v,U))\ (v\cdot x_j | u\cdot x_i,\
 v\cdot x_j\neq u\cdot x_i)\ \ [\ v\cdot x_j\in
 \cup_{u\in U}\,u\,L(u,U)\ ]
$$
the following condition holds:
\begin{equation}
 (\forall w\in \cup_{u\in U}\,u\,L(u,U))
 \ \ [\ u\cdot x_i\not \in wL(w,U\cup \{w\})\ ].
\label{constr}
\end{equation}
} \label{def_constr}
\end{definition}

\begin{definition}
{\em Let $L$ be an involutive division, and $Id(U)$ be a monomial ideal. Then
 an $L-$involutive basis $\bar{U}$ of $Id(U)$ will be called {\em minimal} if
 for any other involutive basis $\bar{V}$ of the same ideal the inclusion
 $\bar{U}\subseteq \bar{V}$ holds.
} \label{def_mmi}
\end{definition}

\begin{proposition}~\cite{GB2}
 If $U\subset \M$ is a finitely generated set with respect to a
 constructive involutive division, then the monomial ideal $Id(U)$ has
 a unique minimal involutive basis.
 \label{pr_mmb}
\end{proposition}

\subsection{Involutive Polynomial Sets}

\begin{definition}
Given a finite set of polynomials $F\subset \R$ and a main ordering
$\succ$, multiplicative and nonmultiplicative variables for $f\in
F$ are defined in terms of $lm(f)$ and the leading monomial set
$lm(F)$.
\end{definition}

The concepts of involutive polynomial reduction and involutive
normal form are introduced similarly to their conventional
analogues~\cite{Buch85,BWK93} with the use of involutive division
instead of the conventional one.

\begin{definition} {\em Let $L$ be an involutive division $L$ on $\M$, and
 let $F$ be a finite set of polynomials. Then we shall say:
 \begin{enumerate}
 \renewcommand{\theenumi}{(\roman{enumi})}
 \item $p$ is {\em $L-$reducible} {\em modulo} $f\in F$ if
  $p$ has a term $t=a\,u$, ($a\in K\setminus \{0\})$, $u\in \M$ such that
  $u=lm(f)\times v$,
  $v\in L(lm(f),lm(F))$. It yields the {\em $L-$reduction} $p\rightarrow
  g=p-(a/lc(f))\,f\, v$.
 \item $p$ is {\em $L-$reducible modulo} $F$ if there is $f\in F$ such
  that $p$ is $L-$reducible modulo $f$.
 \item $p$ is {\em in $L-$normal form modulo $F$} if
  $p$ is not $L-$reducible modulo $F$.
\end{enumerate}
} \label{inv_red}
\end{definition}

\noindent
We denote the $L-$normal form of $p$ modulo $F$ by $NF_L(p,F)$.
In contrast, the conventional normal form will be denoted by $NF(p,F)$.
If monomial $u$ is multiplicative to $lm(f)$ ($f\in F$) and $h=fu$ we shall
write $h=f\times u$.

\begin{definition}{\em A finite polynomial set $F$ is {\em $L-$autoreduced}
 if the leading monomial set $lm(F)$ of $F$ is $L-$autoreduced
 and every $f\in F$ does not contain monomials which are involutively
 multiple of any element in $lm(F)$. }
\end{definition}

\begin{remark} {\em The further definitions and theorems of this section
 which involve the completion ordering $\sqsubset$ generalize those
 in~\cite{GB1} where $\sqsubset$ is the same as the main ordering $\succ$.
 The proofs of the generalized theorems are immediate extensions of the
 underlying proofs in~\cite{GB1}.
}
\end{remark}

\begin{definition}{\em An $L-$autoreduced set $F$ is called
{\em ($L-$)involutive} if
$$
(\forall f\in F)\ (\forall u\in \M)\ \ [\ NF_L(fu,F)=0\ ]\,.
$$
Given $v\in \M$ and an $L-$autoreduced set $F$, if there exist
$f\in F$ such that $lm(f)\sqsubset v$ and
\begin{equation}
(\forall f\in F)\ (\forall u\in \M)\ (lm(f)\cdot u\sqsubset v)\ \
[\ NF_L(fu,F)=0\ ]\,, \label{cond_pinv}
\end{equation}
then $F$ is called {\em partially involutive up to the monomial
$v$} with respect to the ordering $\sqsubset$. $F$ is still said to
be partially involutive up to $v$ if $v\sqsubset lm(f)$ for all
$f\in F$. } \label{def_inv}
\end{definition}

\begin{theorem}\cite{GB1}
 An $L-$autoreduced set $F\subset \R$ is involutive
 with respect to a continuous involutive
 division $L$ iff the following (local) involutivity conditions
 hold
$$
 (\forall f\in F)\ (\forall x_i\in NM_L(lm(f),lm(F)))\
 \ [\ NF_L(f\cdot x_i,F)=0\ ]\,.
$$
Correspondingly, partial involutivity~(\ref{cond_pinv}) holds
iff
$$
(\forall f\in F)\ (\forall x_i\in NM_L(lm(f),lm(F)))\ (lm(f)\cdot
x_i\sqsubset v)\ \ [\ NF_L(f\cdot x_i,F)=0\ ]\,.
$$
\label{th_inv_cond}
\end{theorem}

\begin{theorem}\cite{GB1} If $F\subset \R$ is an $L-$involutive basis of
 $Id(F)$, then
 it is also a \Gr basis, and the equality of the conventional and $L-$normal
 forms $NF(p,F)=NF_L(p,F)$ holds for any polynomial $p\in \R$.
If the set $F$ is partially involutive up to the monomial $v$ with
respect to $\sqsubset$, then the equality of the normal forms
$NF(p,F)=NF_L(p,F)$ holds for any $p$ such that $lm(p)\sqsubset v$.
\label{th_nf}
\end{theorem}

\begin{theorem}
Let $L$ be a continuous involutive division, $F$ be a finite
$L-$autoreduced polynomial set and $NF_L(p,F)$ be an algorithm of
$L-$involutive normal form. Then the following are equivalent:
\begin{enumerate}
 \renewcommand{\theenumi}{(\roman{enumi})}
\item $F$ is an $L-$involutive basis of $Id(F)$.
\item For all $g\in F, x\in NM_L(lm(g),lm(F))$ there is $f\in F$ satisfying
 $lm(g)\cdot x = lm(f)\times w$
 and a chain of polynomials in $F$ of the form
 $$f\equiv f_k,f_{k-1},\ldots,f_0,g_0,\ldots,g_{m-1},g_m\equiv g$$
 such that
 $$ NF_L(S_L(f_{i-1},f_i),F)=NF_L(S(f_0,g_0),F)=NF_L(S_L(g_{j-1},g_j),F)=0\,,
$$
 where $0\leq i\leq k$, $0\leq j\leq m$, $S(f_0,g_0)$ is the conventional
 S-polynomial~\cite{Buch85,BWK93} and $S_L(f_i,f_j)=f_i\cdot x-f_j\times w$
 is its special form which occurs in involutive algorithms~\cite{GB1}.
\end{enumerate}
\label{inv_chain}
\end{theorem}

\noindent
{\bf Proof}\ \ $(i)\Longrightarrow (ii)$ immediately follows from
 Theorems~\ref{th_inv_cond} and~\ref{th_nf} if one takes $f_0=f$, $g_0=g$.
 To prove $(ii)\Longrightarrow (i)$ one suffices to show that
 $NF_L(g\cdot x,F)=0$. Assume for a contradiction that there are
 nonmultiplicative prolongations which are $L-$irreducible to zero modulo $F$.
 Let $g\cdot x$ be such a prolongation which is the lowest with
 respect to the main ordering $\succ$. This means the partial involutivity
 of $F$ up to $lm(g)\cdot x$ with respect to $\succ$. Correspondingly,
 the condition $(ii)$ implies the representation~\cite{BWK93}~(cf. the
 proof of Theorem~8.1. in~\cite{GB1})
 $S_L(f,g)=g\cdot x -f\times w =\sum_{ij}f_iu_{ij}$ where $f_i\in F$
 and $lm(f_iu_{ij})\prec lm(g)\cdot x$ that contradicts
 $NF_L(g\cdot x,F)\neq 0$.
\hfill{\Box}

\begin{corollary}\cite{GB1}
Let $F$ be a finite $L-$autoreduced polynomial set, and let $g\cdot
x$ be a nonmultiplicative prolongation of $g\in F$. If the
following holds
$$
(\forall h\in F)\ (\forall u\in \M)\ \bigl(\ lm(h)\cdot u\sqsubset
lm(g\cdot x)\  \bigr)\ \ [\ NF_L(h\cdot u,F)=0\ ]\,,
$$
$$
 (\exists f,f_0,g_0\in F)
\left[
\begin{array}{l}
 lm(f_0)|lm(f)\,,\ lm(g_0)|lm(g) \\[0.1cm]
 lm(f)|_L lm(g\cdot x)\,,\ lcm(f_0,g_0)\sqsubset lm(g\cdot x) \\[0.1cm]
 NF_L\bigl(f_0\cdot \frac{lt(f)}{lt(f_0)},F\bigl)=
 NF_L\bigl(g_0\cdot \frac{lt(g)}{lt(g_0)},F\bigl)=0
\end{array}
\right]\,,
$$
then the prolongation $g\cdot x$ may be discarded in the course of
an involutive algorithm.
\label{cor_criterion}
\end{corollary}

\begin{remark}{\em Theorem~\ref{th_inv_cond} is the algorithmic
characterization of involutivity whereas Theorem~\ref{th_nf}
relates \Gr bases and involutive bases. Theorem~\ref{inv_chain}
and Corollary~\ref{cor_criterion} yield an involutive analog of the
Buchberger's chain criterion~\cite{Buch85}. }
\end{remark}

\begin{definition}{\em Given a constructive division $L$, a finite
 involutive basis $G$ of ideal $Id(G)$ is called {\em minimal} if
 $lt(G)$ is the minimal involutive basis of the monomial ideal
 generated by $\{lt(f)\ |\ f\in Id(G)\}$.
} \label{min_ipb}
\end{definition}

\begin{theorem}~\cite{GB2}
 A monic minimal involutive basis is unique.
\label{pr_min}
\end{theorem}

\section{Examples of Involutive Divisions}

\subsection{Previously Introduced Divisions}

We give, first, examples of divisions corresponding to separations
introduced by Janet, Thomas and Pommaret for the purpose of
involutivity analysis of PDEs, and two more divisions proposed in
~\cite{GB2}. For the proof of validity of properties (a)-(d) in
Definition~\ref{inv_div} for these divisions we refer
to~\cite{GB1,GB2}.

\begin{definition}Thomas division~{\em \cite{Thomas}.
Given a finite set $U\subset \M$, the variable $x_i$ is considered
as multiplicative for $u\in U$ if $deg_i(u)=max\{deg_i(v)\ |\ v\in
U\}$, and nonmultiplicative, otherwise. } \label{div_T}
\end{definition}

\begin{definition} Janet division~{\em \cite{Janet}. Let the set $U\subset \M$
be finite. For each $1\leq i\leq n$ divide $U$ into groups
labeled by non-negative integers $d_1,\ldots,d_i$:
$$[d_1,\ldots,d_i]=\{\ u\ \in U\ |\ d_j=deg_j(u),\ 1\leq j\leq i\ \}.$$
A variable $x_i$ is multiplicative for $u\in U$ if $i=1$ and
$deg_1(u)=max\{deg_1(v)\ |\ v\in U\}$, or if $i>1$, $u\in
[d_1,\ldots,d_{i-1}]$ and $deg_i(u)=max\{deg_i(v)\ |\ v\in
[d_1,\ldots,d_{i-1}]\}$. } \label{div_J}
\end{definition}

\begin{definition}
 Pommaret division~{\em \cite{Pommaret78}. For a monomial
$u=x_1^{d_1}\cdots x_k^{d_k}$ with $d_k>0$ the variables $x_j,j\geq k$ are
considered as multiplicative and the other variables as nonmultiplicative.
For $u=1$ all the variables are multiplicative.
} \label{div_P}
\end{definition}

\begin{definition} Division I~{\em \cite{GB2}.
Let $U$ be a finite monomial set. The
variable
 $x_i$  is nonmultiplicative for $u\in U$ if there is $v\in U$ such
 that
 $$ x_{i_1}^{d_1}\cdots x_{i_m}^{d_m}u=lcm(u,v),\quad 1\leq m\leq [n/2],
 \quad d_j>0\ \ (1\leq j\leq m)\,,$$
 and $x_i\in \{x_{i_1},\ldots,x_{i_m}\}$.
} \label{div_I}
\end{definition}

\begin{definition} Division II~{\em \cite{GB2}.
For monomial
 $u=x_1^{d_1}\cdots x_k^{d_n}$ the variable $x_i$ is
 multiplicative if $d_i=d_{max}(u)$ where
 $d_{max}(u)=max\{d_1,\ldots,d_n\}$.
} \label{div_II}
\end{definition}

\begin{remark}{\em
 Thomas division, Divisions I and II do not depend on the ordering
 on
 the variables. Janet and Pommaret divisions, as defined, are based on
 the ordering given in~(\ref{var_order}).
 Pommaret division and Division II are globally
 defined in accordance with Definition~\ref{inv_div}.
 }
\end{remark}

\noindent
All these divisions are constructive, and except Pommaret division
they are noetherian \cite{GB1,GB2}.

\subsection{Induced Division}

Now we consider a new class of involutive divisions induced by admissible
monomial orderings~(cf.~\cite{Apel}).

\begin{definition} Induced division. {\em
 Given an admissible monomial
 ordering $\succ $\footnote{This ordering is generally different from
 the main ordering introduced in Sect.2.1.} a variable $x_i$
 is nonmultiplicative for $u\in U$ if there is $v\in U$ such that
 $v\prec u$ and $deg_i(u)<deg_i(v)$.
} \label{ind_div}
\end{definition}

\begin{proposition} The separation given in Definition~\ref{ind_div} is an
 involutive division.
 \label{pr_inddiv}
\end{proposition}

\noindent
{\bf Proof}\ \ Let $L_{\succ}(u,U)$ be the submonoid generated by
 multiplicative variables. We must prove the properties (b-c) in
 Definition~\ref{inv_div} because (a) and (d) hold obviously.

 (b) Let there be a monomial $w$ such that $w\in uL_{\succ}(u,U)\cap
 vL_{\succ}(v,U)$ with $u,v\in U$ and
 $u\neq v$. Assume $u\succ v$ and $\neg v|u$. Then, there is a variable
 $x|(lcm(u,v)/u)$ such that $x\not \in L_{\succ}(u,U)$. Since $v|w$
 we obtain $x|(w/u)$ that contradicts $w\in uL_{\succ}(u,U)$. Thus,
 $v|u$ and $w=v\times (w/v)=v\times [(w/u)(u/v)]$. This yields
 $u\in vL_{\succ}(v,U)$.

 (c) Let $v\in uL_{\succ}(u,U)$ and $w\in vL_{\succ}(v,U)$ with
 $u,v\in U$,
 and, hence, $u|v$ and $u|w$. Suppose $w\not \in uL_{\succ}(u,U)$. It
 follows
 the existence of a variable $x|(w/u)$, $\neg x|(v/u)$ and a monomial
 $t\prec u\prec v$, $t\in U$ such that $x|(lcm(u,t)/u)$. This suggests
 that $x|(lcm(v,t)/v)$ at $t\prec v$, contradicting our initial
 assumption.
\hfill{\Box}

\begin{remark}{\em
 Generally, an ordering
 $\succ $ defining Induced division implies some variable ordering
 which is not compatible with~(\ref{var_order}). However, below we
 assume that the ordering $\succ $ is compatible
 with~(\ref{var_order}).
 }
\end{remark}

\noindent
To distinguish the above divisions, the abbreviations
$T,J,P,I,II,D_{\succ}$ will be sometimes used. For illustrative purposes we
consider three particular orderings to induce involutive divisions:
lexicographical, degree-lexicographical and degree-reverse-lexicographical.
To distinguish these three orderings we shall use the subscripts ${L}$,
${DL}$, ${DRL}$, respectively.

There are certain relations between separations generated by
those divisions.

\begin{proposition}For any $U$, $u\in U$ and
 $\succ$ the inclusions
 $M_T(u,U)\subseteq M_J(u,U)$, $M_T(u,U)\subseteq M_I(u,U)$,
 $M_T(u,U)\subseteq M_{D_{\succ}}(u,U)$ hold. If $U$ is autoreduced
 with respect to Pommaret division, then also
 $M_P(u,U)\subseteq M_J(u,U)$.
 \label{relations}
\end{proposition}

\noindent
{\bf Proof} The inclusion $M_T(u,U)\subseteq
 M_{D_{\succ}}(u,U)$ follows from the observation that $x\in T(u,U)$
 implies $x\in D_{\succ}(u,U)$. The other inclusions
 proved in~\cite{GB1,GB2}.
\hfill{\Box}

\vskip 0.2cm
\noindent
The following example explicitly shows that all eight divisions we
use in this paper are different. In the table we list the
multiplicative variables for every division.

\vskip 0.5cm
\noindent
\begin{example}{\em \cite{GBC98} Multiplicative variables for elements in the
 monomial set
$U=\{x^2y,xz,y^2,yz,z^3\}$ ($x\succ y\succ z$) for different
 divisions:
\begin{center}
\vskip 0.2cm
\bg {|c|c|c|c|c|c|c|c|c|} \hline\hline
Monomial & \multicolumn{8}{c|}{Multiplicative variables}
\\ \cline{2-9}
       & $T$ & $J$  &   $P$   & $I$   & $II$   & $D_{L}$
       & $D_{DL}$ & $D_{DRL}$
         \\ \hline \hline
 $x^2y$ & $x$ & $x,y,z$ & $y,z$ & $x$   & $x$  &  $x$  & $x$ & $x$ \\
 $xz$   & $-$ & $y,z$ &   $z$   & $x$   & $x,z$ &  $x$  & $x,z$ & $x,z$ \\
 $y^2$  & $y$ & $y,z$ & $y,z$   & $y$   & $y$   & $x,y$ & $x,y$ & $y$ \\
 $yz$   & $-$ & $ z $ & $ z $   & $-$   & $y,z$ & $x,y$ & $x,y,z$ & $x,y,z$ \\
 $z^3$  & $z$ & $ z $ & $ z $   & $z$   & $z$   & $x,y,z$ & $z$ & $z$ \\
         \hline \hline
\nd
\vskip 0.5cm
\end{center} \label{exm_1}
}
\end{example}

\begin{proposition} Induced division is noetherian, continuous and
 constructive.
\label{pr_nc_induced}
\end{proposition}

\noindent
{\bf Proof}\ \ {\em Noetherity.} follows immediately from
noetherity of Thomas division and the underlying inclusion in
Proposition~\ref{relations}.

{\em Continuity.} Let $U$ be a finite set, and $\{u_i\}_{(1\leq
i\leq M)}$ be a sequence of elements in $U$ satisfying the
conditions~(\ref{cont_cond}). In accordance with
Definition~\ref{def_cont} we shall show that there are no
coinciding elements in the sequence for each of the two divisions.
There are the following two alternatives:
\begin{equation}
 (i)\ \ u_i=u_{i-1}\cdot x_{j};\qquad (ii)\ \ u_i\neq u_{i-1}\cdot x_{j}\,.
 \label{two_cases}
\end{equation}
Extract from the sequence $\{u_i\}$ the subsequence
$\{t_k\equiv u_{i_k}\}_{(1\leq k\leq K\leq M)}$
of those elements which occur in the left-hand side of
relation $(ii)$ in~(\ref{two_cases}).

Show that $t_{k+1}|_D lcm(t_{k+1},t_k)$ and $\neg t_{k}|t_{k+1}$.
We have $t_{k+1}\times
\tilde{w}_{k+1}=u_{i_k-1}\cdot x_{j_k}=t_{k}\cdot
\tilde{v}_{k}$
where $\neg \tilde{w}_{k+1}|\tilde{v}_{k}$. Indeed, suppose
$\tilde{w}_{k+1} |\tilde{v}_{k}$. Apparently, we obtain the
relation $t_{k+1}=u_{l}\cdot z_l$ where $i_{k}\leq l< i_{k+1}$, and
the variable $x_{j_l}\in NM_D(u_l,U)$, which figures in
Definition~\ref{def_cont} of the sequence $\{u_i\}$, satisfies
$x_{j_l}|\tilde{w}_{k+1}$ and $\neg x_{j_l}|z_l$. This suggests, by
definition of the division, the existence of $p\in U$ such that
$p\prec u_l$ and $deg_{j_l}(u_l) < deg_{j_l}(p)$. Since $p\prec
t_{k+1}$ and $deg_{j_l}(u_l)=deg_{j_l}(t_{k+1})$, it contradicts
multiplicativity of $x_{j_l}$ for $t_{k+1}$.

Therefore, we obtain the relation
$$
\left\{
\begin{array}{l}
t_k\cdot v_k=t_{k+1}\times w_{k+1}\,, \\
gcd(v_{k},w_{k+1})=gcd(v_k,w_k)=1\,,
\end{array}
\right.
$$
which, by Definition~\ref{ind_div}, implies $t_k\succ t_{k+1}$
since $w_{k}\neq 1$ for all $k$.

It remains to prove that elements in the sequence $\{u_i\}_{(1\leq
i\leq M)}$ which occur in the left-hand side of relation $(i)$
in~(\ref{two_cases}) are also distinct. Assume for a contradiction
that there are two elements $u_j=u_k$ with $j < k$. In between
these elements there is, obviously, an element from the left-hand
side of relation $(ii)$ in~(\ref{two_cases}). Let $u_{i_m}$
($j<i_m<k$) be the nearest such element to $u_j$. Considering the
same nonmultiplicative prolongations of $u_k$ as those of $u_j$ in
the initial sequence, one can construct a sequence such that the
subsequence of the left-hand sides of relation (ii)
in~(\ref{two_cases}) has two identical elements $u_{i_k}=u_{i_m}$
with $i_k>i_m$.

{\em Constructivity.} Let $u,u_1\in U$, $v\in D_{\succ}(u_1,U)$ and
 $x_i\not \in D_{\succ}(u,U)$
 be such that $u\cdot x_i=u_1v\times w$, $w\in D_{\succ}(u_1v,U\cup\{u_1v\})$.
 Show that $w\in D_{\succ}(u_1,U)$. Assume that there is
 $x_j\not \in D_{\succ}(u_1,U)$ satisfying $x_j|w$. This implies the
 existence $t\in U$
 satisfying $t\prec u$, $deg_j(t)>deg_j(u_1)$. Then, because $t\prec u_1v$,
 the condition $\neg x_j|v$ leads to the contradictory condition
 $deg_j(t)>deg_j(u_1v)$. Therefore, $x_j|v$.
\hfill{\Box}

\section{Completion of Monomial Sets to Involution}

If $U$ is a finitely generated monomial set with respect to the
fixed involutive division $L$, then its finite completion gives an
involutive basis of the monomial ideal generated by $U$. There may
be different involutively autoreduced bases of the same monomial
ideal. For instance, from Definitions~\ref{div_T} and~\ref{div_J}
it is easy to see that any finite monomial set is Thomas and Janet
autoreduced.  Therefore, enlarging a Thomas or a Janet basis by a
prolongation of any its element and then completing the enlarged
set leads to another Thomas and Janet basis, respectively.
Similarly, Division I and Induced division do not provide
uniqueness of involutively autoreduced bases whereas Pommaret
division and Division II do, as well as any globally defined
division~\cite{GB2}.

\subsection{Completion Algorithm}

\begin{theorem}
 If $U$ is a finitely generated set with respect to a constructive
 involutive division, then the following algorithm computes
 the uniquely defined minimal completion $\bar{U}$ of $U$, that is,
 for any other completion $\tilde{U}$ the inclusion
 $\bar{U}\subseteq \tilde{U}$ holds.
\label{th_compl}
\end{theorem}

\setcounter{cc}{00}
\vskip 0.3cm
\noindent
\hh Algorithm {\bf InvolutiveCompletion:}
\vskip 0.2cm
\noindent
\hh {\bf Input:}  $U$, a finite monomial set
\vskip 0.0cm \noindent
\hh {\bf Output:} $\tilde{U}$, an involutive completion of $U$
\vskip 0.0cm \noindent
\hh \mbegin
\hln
\hhh $\tilde{U}:=U$
\hln
\hhh \mwhile exist $u\in \tilde{U}$ and $x\in NM_L(u,\tilde{U})$ such that
\hln
\hhhh $u\cdot x$ has no involutive divisors in $\tilde{U}$\bb \mdo
\hln
\hhhh \mchoose any $\sqsubset$ and such $u$ and $x$ with the lowest
 $u\cdot x$ w.r.t. $\sqsubset$
\hln
\hhhh $\tilde{U}:=\tilde{U}\cup \{u\cdot x\}$
\hln
\hhh \mend
\hln
\hh \mend
\hln
\vskip 0.3cm

\noindent
{\bf Proof}\ \ This completion algorithm is a slightly generalized
 version of that in~\cite{GB1} where ordering $\sqsubset$ is assumed
 to be fixed in the course of the completion. As proved
 in~\cite{GB1} (see the proof of Theorem 4.14), the output $\bar{U}$ of
 the algorithm and the number of irreducible prolongations are
 invariant on the choice of ordering in line 5.

\begin{corollary} If set $U$ is conventionally
 autoreduced, then the algorithm computes the minimal involutive basis
 of monomial ideal $Id(U)$.
 \label{cor_minmb}
\end{corollary}

\noindent
 In practice, in the course of the completion one has
 to choose the lowest nonmultiplicative prolongation
 and to check whether it has an involutive divisor in the set.
 The timing of computation is thereby determined by the total
 number of prolongations checked.

\begin{theorem} The number of nonmultiplicative prolongations
 checked in the course of algorithm {\bf InvolutiveCompletion} with
 a constructive division $L$ is invariant on the choice of completion
 ordering in line 5.
\label{th_invariant}
\end{theorem}

\vskip 0.2cm
\noindent
{\bf Proof}\ \ As it has been noticed in the proof of
Theorem~\ref{th_compl}, the number of irreducible prolongations, as
well as the completed set itself, is invariant. Therefore, we must
prove invariance of the number of reducible prolongations.

Let there be two different completion procedures of $U$ to
$\bar{U}$ based on different choice of completion orderings. Assume
that the first procedure needs more reducible prolongations to
check than the second one. Let $u\cdot x=v\times w$ $(u,v\in
\tilde{U}_1)$ be the first prolongation checked in the course of the
first procedure and such that in the course of the second one the
prolongation is not checked. This suggests $x\in
M_L(u,\tilde{U}_2)$ where $\tilde{U}_2$ is the current set for the
second procedure. If $w\neq 1$, by admissibility of a completion
ordering, we obtain $u\times x=v\times w$ $(u,v\in
\tilde{U}_2)$. From property (b) in Definition~\ref{inv_div} we
deduce $u\in vL(v,\tilde{U}_2)$, and, hence, $x$ cannot be
nonmultiplicative for $u$ as we assumed for the first procedure.

If $w=1$ we find that $x\in NM_L(u,\tilde{U}_1)\cap
M_L(u,\tilde{U}_2)$ where $u$ and $v=ux$ are elements in both
$\tilde{U}_1,\tilde{U}_2$. From the property (d) in
Definition~\ref{inv_div} and invariance of the final completed set
$\bar{U}$ is follows that in some step of the second procedure $x$
becomes nonmultiplicative for $u$. Then the prolongation $u\cdot x$
will be also checked that contradicts our assumption.

\hfill{\Box}

\vskip 0.2cm
\noindent
This theorem generalizes Remark~3.13 in~\cite{GBC98} which is
concerned with $L-$autoreduced sets and fixed completion orderings.

\begin{example} {\em (Continuation of Example~\ref{exm_1}). The
 minimal involutive bases of the ideal generated by
 the set $U=\{x^2y,xz,y^2,yz,z^3\}$ ($x\succ y\succ z$) are given by
\begin{eqnarray*}
\bar{U}_T&=&\{x^2y^2z^3,x^2y^2z^2,x^2y^2z,x^2y^2,x^2yz^3,
   x^2yz^2,x^2yz,x^2y,x^2z^3, \\
&&   x^2z^2,x^2z,xy^2z^3,xy^2z^2,
   xy^2z,xy^2,xyz^3,xyz^2,xyz,xz^3,xz^2, \\
&&   xz,y^2z^3, y^2z^2,y^2z,y^2,yz^3, yz^2,yz,z^3\}\,,\\
\bar{U}_J&=&\{x^2y, x^2z, xy^2, xyz, xz, y^2, yz, z^3\}\,,\\
\bar{U}_P&=&\{x^2y,x^2z,xy^2,xyz,xz,y^2,yz,z^3,\ldots,x^ky,\ldots,x^lz,\ldots \}\,, \\
\bar{U}_I&=&\{x^2y^2z^3, x^2y^2z^2, x^2y^2z, x^2y^2, x^2yz^3, x^2yz^2,
          x^2yz, x^2y, xy^2z^3, \\
&&        xy^2z^2, xy^2z, xy^2,
     xyz^3, xyz^2, xyz, xz^3, xz^2, xz, y^2z^3, y^2z^2, \\
&&     y^2z, y^2, yz^3, yz^2, yz, z^3\}\,,\\
\bar{U}_{II}&=&\{x^2y^2, x^2y, xy^2, xyz, xz, y^2, yz, z^3\}\,, \\
\bar{U}_{L}&=&\{x^2y, xz^2, xz, y^2, yz^2, yz, z^3\}\,, \\
\bar{U}_{DL}&=&\{x^2y, xz, y^2, yz, z^3\}\,, \\
\bar{U}_{DRL}&=&\{x^2y, xy^2, xz, y^2, yz, z^3\}\,,
\end{eqnarray*}
where $k,l\in \N$ ($k,l>2$), and subscripts in the left-hand sides
stand for different involutive divisions considered in Section 3.
This example explicitly shows that Pommaret division is not noetherian,
since it leads to an infinite monomial basis.
} \label{exm_2}
\end{example}

\subsection{Pair Property}

In the course of algorithm {\bf InvolutiveCompletion} the current
monomial set $\tilde{U}$ is enlarged by irreducible
nonmultiplicative prolongations in line 6. As this takes place, for
a nonglobally defined division one has to recompute the separation
into multiplicative and nonmultiplicative variables for {\em all}
monomials. The next definition and proposition give a prescription
for efficient recomputing.

\begin{definition}
{\em We shall say that an involutive division $L$ is {\em pairwise}
if for any finite set $U$ and any $u\in U$ $(U\setminus \{u\}\neq
\emptyset)$, the following holds:
$$
 L(u,U)=\cap_{v\in U\setminus \{u\}} L(u,\{v\})
$$
or, equivalently,
\begin{equation}
 M_L(u,U)=\cap_{v\in U\setminus \{u\}} M_L(u,\{v\}),\qquad
 NM_L(u,U)=\cup_{v\in U\setminus \{u\}}NM_L(u,\{v\}).
 \label{pairwise}
\end{equation}
\label{def_pairwise}
}
\end{definition}

\noindent
Therefore, for a pairwise division $L$ and a monomial set $U$ the
correction of the separation due to enlargement of $U$ by an
element $v$ is performed by formula
\begin{equation}
NM_L(u,U\cup \{v\})=NM_L(u,U)\cup NM_L(u,\{u,v\}).
\label{recomputing}
\end{equation}

\begin{proposition}
 All the above defined divisions are pairwise.
 \label{pr_pairwise}
\end{proposition}

\noindent
{\bf Proof}\ \ Pommaret division and Division II, as globally
 defined divisions, are trivially pairwise.

 {\em Thomas division}. Since
 $$max \{deg_i(w)\ |\ w\in U\}=max_{v\in U\setminus \{u\}}\{max
 \{deg_i(u),deg_i(v)\}\},$$
 Definition~\ref{div_T} implies apparently~(\ref{pairwise}).

 {\em Janet division}. For $i=1$, by our convention~(\ref{var_order}) and
  Definition~\ref{div_J}, Janet case is reduced to Thomas
 one, and we are done. Let now $i>1$ and
 $deg_1(u)=d_1,\ldots,deg_{i-1}(u)=d_{i-1}$. If the group
 $[d_1,\ldots,d_{i-1}]$ of elements in $U$ contains, in addition to $u$,
 some extra elements, then $$max \{deg_i(w)\ |\ w\in [d_1,\ldots,d_{i-1}]\}=
 max_{v\in [d_1,\ldots,d_{i-1}]\setminus \{u\}}\{max
 \{deg_i(u),deg_i(v)\}\},$$ and $x_i\in M_J(u,U)$, otherwise. This
 suggests the pairwise property.

 {\em Division I and Induced division}. For these division the pairwise
 property follows immediately from Definitions~\ref{div_I}
 and~\ref{ind_div}.
\hfill{\Box}

\subsection{Monotonicity}

Consider now another optimization related to the choice of a
nonmultiplicative prolongation in line 5 of algorithm {\bf
InvolutiveCompletion}. The choice of the lowest prolongation with
respect to some fixed ordering $\sqsubset$ is called {\em normal
selection strategy}~\cite{GB1}.

\begin{definition}{\em Given a division $L$ and an admissible ordering
$\sqsubset$, a monomial set $U$ will be called {\em complete up to
monomial $w$ with respect to $\sqsubset$} if
\begin{equation}
(\forall u\in U)\ (\forall x\in NM_L(u,U))\ (u\cdot x\sqsubseteq
w)\
\ [\ u\cdot x\in C_L(U)\ ]\,,
\label{part_compl}
\end{equation}
where $C_L(U)$ is involutive cone of $U$ by Definition~\ref{cone}.
We call monomial $w$ {\em bound of completeness} for $U$. If
$u\cdot x\sqsupset w$ for all $u\in U$, $x\in NM_L(u,U)$, then we
shall still say that $U$ is completed up to $w$.}
\label{def_partcompl}
\end{definition}

\begin{definition}{\em We shall call division $L$ {\em monotone for
$\sqsubset$} if for any set $U$ and any monomial $w\in \M$
satisfying~(\ref{part_compl}) the following holds:
$$
 (\forall v\in U)\ (\forall x\in NM_L(v,U))\ (v\cdot x\not \in C_L(U))
\ \ [\ U\cup \{v\cdot x\}\ {is\ complete\ up\ to\ } w\ ].
$$
We shall say that $L$ is {\em monotone} if its monotonicity holds
for any ordering $\sqsubset$.
 } \label{def_mon}
\end{definition}
Thus, monotonicity means that enlargement of $U$ by an irreducible
nonmultiplicative prolongation does not decrease its completeness
bound.

\begin{remark}{\em If a division $L$ is monotone for an ordering
$\sqsubset $, then the choice of the latter as a completion
ordering is beneficial for the algorithm {\bf
InvolutiveCompletion}. By Theorem~(\ref{th_invariant}), the total
number of prolongations checked is invariant on the ordering.
Monotonicity of the latter allows one to omit recomputing
separations and checking prolongations which are lower than the
current completeness bound.}
\label{rem_opt}
\end{remark}

\noindent
Now we consider the monotonicity properties of different divisions
defined in Sect.2. Pommaret division and Division II, as globally
defined, are trivially monotone.

\begin{proposition} Thomas division is monotone.
\label{mon_T}
\end{proposition}

\noindent
{\bf Proof}\ \ From Definition~\ref{div_T} it follows immediately
that $T(u,U)=T(u,U\cup \{v\cdot x\})$ for any $v\in U$, $x\in
 NM_T(v,U)$.
\hfill{\Box}

\begin{proposition} Janet division is monotone for lexicographical
ordering.
\label{mon_J}
\end{proposition}

\noindent
{\bf Proof}\ \ Denote the lexicographical completion ordering
compatible with~(\ref{var_order}) by $\sqsubset_{Lex}$. Consider a
nonmultiplicative prolongation $v\cdot x_j\not \in C_J(U)$ ($v\in
U$) such that $v\cdot x_j\sqsupset_{Lex} w$ where $w$ is the
completeness bound of $U$ in accordance with~(\ref{part_compl}).

Suppose there is a pair $\{u\in U,x_k\}$, satisfying
\begin{equation}
ux_k\in C_J(U),\quad u\cdot x_k\not \in C_J(U\cup
\{vx_j\}),\quad u\cdot x_k\sqsubseteq_{Lex} w\sqsubset_{Lex}
v\cdot x_j,
\label{relat}
\end{equation}
and consider the lowest such pair with respect to
$\sqsubset_{Lex}$.

If $x_k\in J(u,U)$ we obtain
$$
\left\{
\begin{array}{ll}
deg_1(u)=deg_1(v)+1,\quad deg_k(u) < deg_k(v) & if\ j=1, \\
deg_i(u)=deg_i(v)\ (i<j),\quad deg_j(u)=deg_j(v)+1,\quad deg_k(u) <
deg_k(v) & if\ j>1.
\end{array}
\right.
$$
Here $k>j$ and, if $k-j>1$, then $deg_m(u)=deg_m(v)$ for all
$l<m<k$. Consider now two alternatives:

(i) $w\in U$. In this case conditions~(\ref{relat}) are
contradictory since from the rightmost condition it follows
$deg_i(u)=deg_i(w)$ $(i<k)$ and $deg_k(w)>deg_k(u)$, that is,
$x_k\in NM_J(u,U)$.

(ii) $w\not \in U$. Then there is $t\in U$ such that $w\in
tJ(t,U)$. Because $x_k\in J(u,U)$, for some $1\leq p<k$ we have
$deg_i(t)=deg_i(u)=deg_i(vx_j)$ where $i\leq p$ and
$deg_{p+1}(t)<deg_{p+1}(u)\leq deg_{p+1}(w)$. Thus we obtain
contradiction with $x_{p+1}\in J(t,U)$ which follows from $w\in
tJ(t,U)$.

It is remains to prove that if $x_k\in NM_J(u,U)$, then $u\cdot
x_k\in C_J(U\cup \{v\cdot x_j\})$. If $u\cdot x_k\in U$ we are
done. Otherwise, we have $u\cdot x_k=q_1\,r_1$ for some $q_1\in U$,
$r_1\in J(q,U)$ and $r_1\not \in J(U\cup \{v\cdot x_j\})$. Hence,
there is $x_{i_1}|r_1$, $x_{i_1}\in NM_J(q_1,U\cup \{v\cdot
x_j\})$, and $deg(q_1\cdot x_{i_1})\sqsubset_{Lex}deg(u\cdot x_k)$.
Then, by our assumption that prolongation $u\cdot x_k$ is the
lowest satisfying ~(\ref{relat}), we have $q_1\cdot
x_{i_1}=q_2\times r_2$, $q_2\in U$, $r_2\in J(q_2,U\cup v\cdot
x_j)$. By property (d) in Definition~\ref{inv_div}, it yields
$r_2\in J(q_2,U)$, and, hence, $q_1\times x_{i_1}=q_2\times r_2$ in
$U$. This is impossible, because any monomial set is Janet
autoreduced.
\hfill{\Box}

\begin{remark}{\em Janet division is not monotone for
 degree-lexicographical and degree-reverse-lexicographical
 orderings as the following example shows.
} \label{rem_nm_J}
\end{remark}

\begin{example}{\em Consider the conventionally autoreduced
 set $U=\{xz^2,x^2z,yzt^2\}$. Let completion ordering
 $\sqsubset $ be degree-lexicographical or degree-reverse-lexicographical
 ordering with $x\sqsupset y\sqsupset z\sqsupset t$. $U$ is complete up to
 $w=x^2yz$. The lowest irreducible prolongation is $xyzt^2\sqsupset w$.
 The next one in the set $\{xz^2,x^2z,yzt^2,xyzt^2\}$ is
 $xyz^2\sqsubset w$.
}
\label{ex_counterJanet}
\end{example}

\begin{example}{\em Consider the set $U=\{xy^2w^2,xzt,yzt\}$ and Division I
generating the separation:
\begin{center}

\vskip 0.2cm
\bg {|l|c|c|} \hline\hline
Monomial & \multicolumn{2}{c|}{Division I}
 \\ \cline{2-3}
           & $M_I$ & $NM_I$  \\ \hline
 $p=xy^2w^2$ & $x,y,w$ & $z,t$  \\
 $v=xzt$  & $x,z,t$  & $y,w$\\
 $u=yzt$   & $y,z,t,w$ & $x$ \\ \hline \hline
\nd
\vskip 0.2cm
\end{center}

\noindent
Let $x\sqsupset y\sqsupset z\sqsupset t\sqsupset w$ and $\sqsupset$
be any of the orderings: lexicographical, degree-lexicographical or
degree-reverse-lexicographical. The bound of completeness for $U$
is $xyw^2$. We find that $v\cdot w=u\cdot x=xztw$ is the lowest
irreducible prolongation in $U$, and the next one for $U\cup
\{xztw\}$ is $u\cdot w\sqsubset xyw^2$. Therefore, Division I is not
monotone for three orderings considered.
}
\label{ex_counterDivI}
\end{example}

\begin{proposition} Induced division is monotone for the ordering
 which induces this division.
\label{pr_monind}
\end{proposition}

\noindent
{\bf Proof}\ \ By Definition~\ref{ind_div} of $D_{\sqsubset}$,
enlargement of $U$ by irreducible nonmultiplicative prolongation
$v\cdot x_j\sqsupset w$ $(v\in U)$ does not change the reducibility
properties of those prolongations $ux_k$ $(u\in U)$ which satisfy
$ux_k\sqsubset v\cdot x_j$.
\hfill{\Box}

\section{Construction of Involutive Bases for Polynomial Ideals}

In this section we present the following algorithm for computation
of minimal involutive bases of polynomial ideals which generalizes the
algorithm of paper~\cite{GB2} to different completion and main orderings.

\vskip 0.3cm
\setcounter{cc}{00}
\noindent
\h Algorithm {\bf MinimalInvolutiveBasis:}
\vskip 0.2cm
\noindent
\h {\bf Input:} \begin{minipage}[t]{12cm}
  $F$, a finite polynomial set; $L$, an involutive division;\\
  $\succ$, a main ordering; $\sqsubset$, a completion ordering
\end{minipage}
\vskip 0.1cm
\noindent
\h {\bf Output:} $G$, the minimal involutive basis of $Id(F)$ if algorithm
terminates
\vskip 0.0cm
\noindent
\h \mbegin
\hln
\hh $F:=Autoreduce(F)$
\hln
\hh \mchoose $g\in F$ with the lowest $lm(g)$ w.r.t. $\prec$
\hln
\hh $T:=\{(g,lm(g),\emptyset)\}$;\ \ $Q:=\emptyset$;\ \ $G:=\{g\}$
\hln
\hh \mfore $f\in F\setminus \{g\}$\bb \mdo
\hln
\hhh $Q:=Q\cup \{(f,lm(f),\emptyset )\}$
\hln
\hh \mrepeat
\hln
\hhh $h:=0$
\hln
\hhh \mwhile $Q\neq \emptyset$\bb \mand $h=0$\bb \mdo
\hln
\hhhh \mchoose $g$ in $(g,u,P)\in Q$ with the lowest $lm(g)$ w.r.t. $\prec$
\hln
\hhhh $Q:=Q\setminus \{(g,u,P)\}$
\hln
\hhhh \mif $Criterion(g,u,T)$ is false \bb \mthen $h:=NF_L(g,G)$
\hln
\hhh \mif $h\neq 0$\bb \mthen $G:=G\cup \{ h \}$
\hln
\hhhh \mif $lm(h)=lm(g)$\bb \mthen $T:=T\cup \{(h,u,P\cap NM_L(h,G))\}$
\hln
\hhhh \melse $T:=T\cup \{(h,lm(h),\emptyset )\}$
\hln
\hhhhh \mfore  $(f,v,D)\in T$ s.t. $lm(f) \succ lm(h)$\bb \mdo
\hln
\hhhhhh $T:=T\setminus \{(f,v,D)\}$;\ \ $Q:=Q \cup \{(f,v,D)\}$;\ \
    $G:=G\setminus \{f\}$
\hln
\hhhhh \mfore  $(f,v,D)\in T$ \bb \mdo
\hln
\hhhhhh $T:=T\setminus \{(f,v,D)\}\cup  \{(f,v,D\cap NM_L(f,G))\}$
\hln
\hhh \mwhile exist $(g,u,P)\in T$ and $x\in NM_L(g,G)\setminus P$ and,
     if $Q \neq \emptyset$,
\hln
\hhhh s.t. $lm(g\cdot x) \prec lm(f)$ for all $f$ in $(f,v,D)\in Q$\bb
\mdo
\hln
\hhhh \mchoose such $(g,u,P),x$ with the lowest $lm(g)\cdot x$ w.r.t. $\sqsubset$
\hln
\hhhh $T:=T\setminus \{(g,u,P)\} \cup \{(g,u,P\cup \{x\})\}$
\hln
\hhhh \mif $Criterion(g\cdot x,u,T)$ is false \bb \mthen $h:=NF_L(g\cdot x,G)$
\hln
\hhhhh \mif $h\neq 0$\bb \mthen $G:=G\cup \{h\}$
\hln
\hhhhhh \mif $lm(h)=lm(g\cdot x)$ \bb \mthen $T:=T\cup \{(h,u,\emptyset )\}$
\hln
\hhhhhh \melse $T:=T\cup \{(h,lm(h),\emptyset )\}$
\hln
\hhhhhhh \mfore  $f$ in $(f,v,D)\in T$\bb with\bb $lm(f) \succ lm(h)$\bb \mdo
\hln
\hhhhhhhh $T:=T\setminus \{(f,v,D)\}$;\ \ $Q:=Q \cup \{(f,v,D\})$;\ \
 $G:=G\setminus \{f\}$
\hln
\hhhhhhh \mfore  $(f,v,D)\in T$ \bb \mdo
\hln
\hhhhhhhh $T:=T\setminus \{(f,v,D)\}\cup  \{(f,v,D\cap NM_L(f,G))\}$
\hln
\hh \muntil $Q\neq \emptyset$
\hln
\h \mend
\hln
\vskip 0.2cm
\noindent
$Criterion(g,u,T)$ is true provided that if there is $(f,v,D)\in T$
such that $lm(f)|_Llm(g)$ and $lcm(u,v) \sqsubset lm(g)$. Correctness
of this criterion, which is just the involutive form~\cite{GB1} of
Buchberger's chain criterion~\cite{Buch85}, is provided by
Corollary~\ref{cor_criterion}.

\begin{theorem}Let $F$ be a finite subset of $\R$ and $L$ be a
 constructive involutive division. Suppose the main ordering
 $\succ$ is degree compatible. Then the algorithm {\bf
 MinimalInvolutiveBasis} computes a minimal involutive basis of $Id(F)$
 if this basis is finite. If $L$ is noetherian, then the basis
 is computed for any main ordering.
\label{th_alg_min}
\end{theorem}

\noindent
{\bf Proof}\ \ The proof is the same as in~\cite{GB2} and based on
Theorems~\ref{th_inv_cond}, \ref{th_nf} and \ref{th_compl},
Corollaries \ref{cor_criterion} and \ref{cor_minmb}.
\hfill{\Box}

\begin{proposition}
The conventional autoreduction of the input polynomial set in line
2 is optional and may be omitted.
\label{autoreduction}
\end{proposition}

\noindent
{\bf Proof}\ \ Let $F$ be a non-autoreduced set and the algorithm
start with line 3. Subsequent to the initialization in lines 4-6
the upper {\bf while}-loop selects, first of all, those polynomials
in the triple set $Q$ which have the same leading term as the
element in $G=\{g\}$. If there is such a polynomial in the triple
set $Q$ with nonzero involutive normal form $h$ computed in line
12, then $lm(h)\prec lm(g)$. It follows from lines 13 and 17
that $G$ becomes the one-element set $\{h\}$ as an input for the
lower {\bf while}-loop.

Thus, by restriction in line 21 for nonmultiplicative prolongations
checked and redistribution of polynomials in line 29, in every step
of the algorithm we have $lm(g)\prec lm(f)$ for any $g$ in
$(g,u,P)\in T$ and $f$ in $(f,v,D)\in Q$ whenever the set $Q$ is
nonempty.

Furthermore, as proved in~\cite{GB1,GB2}, in some step of the
algorithm a polynomial $h$ is added to the current polynomial set
$G$ in line 13 or in line 25, such that $h$ is an element in the
reduced \Gr basis of $Id(F)$ with the lowest leading monomial with
respect to the main ordering $\prec$. It implies the
reduction of $G$ to the one-element set $G=\{h\}$, and transfer of
the rest to $Q$. Then $G$ is sequentially completed by other
polynomials from the reduced \Gr basis and their nonmultiplicative
prolongations. In so doing, the completion of $lm(G)$ due to the
redistribution of polynomials between sets $T$ and $Q$ in lines 17
and 29 is monotone with respect to $\prec$.

Therefore, the output of algorithm {\bf MinimalInvolutiveBasis}
irrespective of autoreduction in line 2 is the same as it would be
for the reduced \Gr basis in the input.
\hfill{\Box}

\begin{remark} {\em The choice of a completion ordering which is
monotone for $L$ preserves, obviously, the partial involutivity of
the intermediate polynomial set $G$ in the case of its enlargement
in line 25, if $lm(h)=lm(g\cdot x)$. Therefore, similar to the
monomial case (c.f. Remark~\ref{rem_opt}), this saves computing
time for recomputing separations and checking irreducibility of
nonmultiplicative prolongations unless $L$ is globally defined
anyway. }
\label{rem_opt_pol}
\end{remark}

\section{Conclusion}

The above described optimizations concern only that part of
computing involutive bases which is related to completion by
nonmultiplicative prolongations with irreducible leading terms.
Another important step is to search for an involutive divisor among
the leading monomials of an intermediate basis. This is important
for efficient computation of the involutive normal form in lines 12
and 24 of algorithm {\bf MinimalInvolutiveBasis}. Some related
optimizations are considered in~\cite{GBC98} for the purpose of
implementing the algorithm {\bf InvolutiveCompletion} in
Mathematica for divisions of Sect.3.

A promising way to the further optimization of computation is
related to the ideas of paper~\cite{Apel}. By appropriate dynamical
refinement of an involutive division in the course of computation,
one can decrease the total number of nonmultiplicative
prolongations to be checked. This may lead to a notable reduction
of computing time.

Algorithm {\bf MinimalInvolutiveBasis} has been implemented in
Reduce for Pommaret division. Computer experiments showed that this
algorithm is somewhat faster than our previous version of
involutive algorithm also implemented in Reduce for Pommaret
bases~\cite{GB1}. For a nonglobally defined division the difference
in speed is to be much greater as algorithm {\bf
MinimalInvolutiveBasis} deals with fewer intermediate polynomials
and avoids intermediate autoreductions~\cite{GB2}.

With the new implementation one needs, for example, 57 seconds to
compute a degree-reverse-lexicographical Pommaret basis for 6th
cyclic roots on an Pentium 100 Mhz computer, and 30 seconds for the
6th Katsura system. By comparison, the PoSSo software for computing
\Gr bases\footnote{URL: http://janet.dm.unipi.it/posso\_demo.html}
needs for these examples 24 and 36 seconds, respectively.

\section{Acknowledgements} The author thanks Joachim Apel,
Matthias Berth and Yuri Blinkov for helpful discussions and
comments. This work was partially supported by grant INTAS-96-0842
and grants from the Russian Foundation for Basic Research
No. 98-01-00101, 00-15-96691.

\end{document}